\documentclass[11 pt]{amsart}
\usepackage[applemac]{inputenc}

\usepackage[dvipsnames,svgnames,x11names]{xcolor}
\definecolor{ffmblue}{HTML}{006092}
\usepackage[top=20mm,bottom=20mm,left=20mm,right=20mm]{geometry}

\usepackage[hyphens]{url}
\usepackage[colorlinks=true,linkcolor=Maroon,citecolor=Maroon,urlcolor=ffmblue,hypertexnames=false,linktocpage]{hyperref}
\usepackage{bookmark}
\usepackage{amsmath,thmtools,mathtools}
\mathtoolsset{showonlyrefs=true}
\newcounter{mparcnt}

\usepackage{fancyhdr}
\usepackage{esint}
\usepackage{enumerate}

\usepackage{pictexwd,dcpic}
\usepackage{graphicx}
\usepackage{caption}
\swapnumbers
\declaretheorem[name=Theorem,numberwithin=section]{thm}
\declaretheorem[name=Remark,style=remark,sibling=thm]{rem}
\declaretheorem[name=Lemma,sibling=thm]{lemma}
\declaretheorem[name=Proposition,sibling=thm]{prop}
\declaretheorem[name=Definition,style=definition,sibling=thm]{defn}
\declaretheorem[name=Corollary,sibling=thm]{cor}
\declaretheorem[name=Assumption,style=definition,sibling=thm]{assum}

\declaretheorem[name=Notation,style=definition,sibling=thm]{nota}

\numberwithin{equation}{section}

\newcommand{\sub}{\subset}

\newcommand{\bbR}{\mathbb{R}}

\newcommand{\bbS}{\mathbb{S}}

\newcommand{\8}{\infty}

\newcommand{\de}{\delta}

\newcommand{\ka}{\kappa}

\newcommand{\si}{\sigma}

\newcommand{\Ga}{\Gamma}

\newcommand{\del}{\partial}
\newcommand{\n}{\nabla}

\newcommand{\ip}[2]{\left\langle #1,#2 \right\rangle}
\newcommand{\fr}[2]{\frac{#1}{#2}}
\newcommand{\tfr}[2]{\tfrac{#1}{#2}}

\DeclareMathOperator{\dist}{dist}

\DeclareMathOperator{\vol}{vol}

\newcommand{\pf}[1]{\begin{proof}#1 \end{proof}}
\newcommand{\eq}[1]{\begin{equation}\begin{alignedat}{2} #1 \end{alignedat}\end{equation}}

\newcommand{\br}[1]{\left(#1\right)}
\newcommand{\abs}[1]{\lvert #1\rvert}
\newcommand{\enum}[1]{\begin{enumerate}[(i)] #1 \end{enumerate}}

\newcommand{\mt}{\mapsto}

\newcommand{\hp}{\hphantom}
\newcommand{\q}{\quad}

\begin{document}
\title[Quermassintegral inequalities in the sphere]{The quermassintegral inequalities for horo-convex domains in the sphere}
\begin{abstract}
We study a new notion of convexity for subsets of the unit sphere, which closely resembles the horo-convexity for subsets of the hyperbolic space. We call this notion, accordingly, horo-convexity. For horo-convex hypersurfaces of the unit sphere, we prove the smooth convergence of the classical Guan/Li flow of inverse type and use this result to prove the full set of quermassintegral inequalities for horo-convex hypersurfaces of the unit sphere. 
\end{abstract}
\date{\today}
\author{Shujing Pan}
\address{\flushleft\parbox{\linewidth}{{\bf Shujing Pan}\\University of Science and Technology of China\\ School of Mathematical Sciences\\ 230026 Hefei \\ P.R. China\\ \\Goethe-Universit\"at\\ Institut f\"ur Mathematik\\ Robert-Mayer-Str.~10\\ 60325 Frankfurt\\ Germany\\ {\href{mailto:pan@math.uni-frankfurt.de}{pan@math.uni-frankfurt.de}}}}

\author{Julian Scheuer}
\address{\flushleft\parbox{\linewidth}{{\bf Julian Scheuer}\\Goethe-Universit\"at\\ Institut f\"ur Mathematik\\ Robert-Mayer-Str.~10\\ 60325 Frankfurt\\ Germany\\ {\href{mailto:scheuer@math.uni-frankfurt.de}{scheuer@math.uni-frankfurt.de}}}}

\maketitle

\section{Introduction}

\subsection*{Results}

The setting of this paper is the open northern hemisphere $\bbS^{n+1}_{+}$ of the standard round sphere of radius $1$, equipped with the round metric $\bar g = \ip{\cdot}{\cdot}$.
This paper deals with inverse-type locally constrained curvature flows in $\bbS^{n+1}_{+}$ given by
\eq{\label{flow}\del_{t}x=\br{\fr{\phi'(r)}{F}-u}\nu,}
where $\phi = \sin$ and $r$ is the spherical distance to the north pole.
 The function 
\eq{u = \ip{\phi\del_{r}}{\nu}}
is called support function and $F$ is a function depending on the principal curvatures of the flowing surfaces and satisfying certain properties to be specified below. Of particular interest are the quotients $H_{k}/H_{k-1}$ of normalised elementary symmetric polynomials, because with this choice the monotonicity properties of the flow work in favour towards a proof of the quermassintegral inequalities. For convex hypersurfaces of the unit sphere, the full set of these inequalities is still open, and a nontrivial curvature condition, which allows a proof of the full set of these inequalities, is  missing. In this paper, we solve the latter question by defining a new natural convexity condition, called horo-convexity, which geometrically resembles the notion of horo-convexity in the hyperbolic space. Hence this paper provides a natural spherical counterpart to the paper \cite{WangXia:07/2014}, which proves the quermassintegral inequalities for horo-convex hypersurfaces in the hyperbolic space. Before we go deeper into the history, we formulate our results.

The new kind of convexity we are using resembles the horo-convexity in the hyperbolic space in the sense that it can be geometrically defined by the existence of outscribed balls which lie entirely in the upper hemisphere. In the hyperbolic ball model this would precisely give the horo-convexity. In the following, we give the analytical definition, where $h$ is the second fundamental form and $g$ the induced metric on the hypersurface.
Our convention for $h$ comes from the following Gaussian formula
\eq{\bar\n_{X}Y = \n_{X}Y - h(X,Y)\nu}
for vector fields $X$ and $Y$ along a strictly convex hypersurface, and where $\nu$ is the normal pointing out of the domain enclosed by that hypersurface.

\begin{defn}
A $C^{2}$-hypersurface of $(\bbS^{n+1}_{+},\bar g)$ is called {\it horo-convex}, if it satisfies
\eq{\phi' h\geq (1-u)g.}
\end{defn}

Next we collect our assumptions on the curvature function $F$. 
\begin{assum}\label{assum:F}
We assume $F\in C^{\8}(\Ga_{+})$, where $\Ga_{+}\sub \bbR^{n}$ is the positive cone, satisfies the following assumptions:
\enum{
\item $F$ is strictly monotone increasing in each argument;
\item $F(1,\dots,1)=1$;
\item $F$ is homogeneous of degree $1$;
\item $F$ is concave;
\item $F$ is inverse concave, i.e. the function $\ka = (\ka_{i})\mt 1/F(\ka_{i}^{-1})$ is concave.
}
\end{assum}

Our first result is the smooth convergence of the standard inverse-type Guan/Li flow, given the initial horo-convexity. Also this result is entirely new and interesting in its own right, because unless for the special case $F = H_{n}/H_{n-1}$, no convergence results are known until today. Remarkably, notice that we do not assume the initial hypersurface to be starshaped around the north pole, which used to be a standard assumption in previous works on the topic. This underpins the naturality of our convexity condition.

\begin{thm}\label{thm:flow}
Let $M_{0}\sub \bbS^{n+1}_{+}$ be horo-convex. Let $F$ satisfy \autoref{assum:F}. Then the flow \eqref{flow} starting from $M_{0}$ exists for all times and converges smoothly to a sphere centred at the north pole.
\end{thm}

Of course, the main application is the full set of quermassintegral inequalities for horo-convex hypersurfaces of the sphere.

\begin{thm}\label{thm:QM}
Let $M\sub \bbS^{n+1}_{+}$ be horo-convex. Then there holds
\eq{\label{thm:QM-A}W_{k}(M)\geq f_{k}(W_{k-1}(M)),\q 1\leq k\leq n,}
where $f_{k}$ is a function which gives equality on geodesic spheres.
The equality case is precisely achieved on geodesic spheres.
\end{thm}

The functionals $W_{k}$ are the well-known quermassintegrals for convex sets in the sphere and they are defined by
\eq{W_{0}(M) = \vol(\hat M),\q W_{1}(M) = \tfr{1}{n+1}\abs{M}}
and
\eq{\label{Wk}W_{k+1}(M) = \fr{1}{n+1}\int_{M}H_{k} + \fr{k}{n+2-k}W_{k-1}(M),}
where $\hat M$ is the domain enclosed by $M$ and $\abs{M}$ is the surface area of $M$.
We do not give another introduction to these functionals. They have been discussed a lot and further details can be found in \cite{Solanes:/2006}.

\begin{proof}[Proof of \autoref{thm:QM} using \autoref{thm:flow}]
It is well known, that under a variation of the form $\del_{t}x = \Phi\nu$, the quermassintegrals evolve like
\eq{\del_{t}W_{k}(M_{t}) = c_{n,k}\int_{M_{t}}\Phi H_{k},}
where $c_{n,k}$ are suitable constants, cf. \cite[Prop.~2.9]{ChenGuanLiScheuer:/2022}. A straightforward computation, which employ the Hsiung-Minkowski integral identities \cite[Prop.~2.5]{GuanLi:/2015},
\eq{\int_{M_{t}}\phi'H_{k-1} =\int_{M_{t}}uH_{k},}
and the Newton-Maclaurin inequalities
\eq{H_{k-1}H_{k+1}\leq H_{k}^{2}}
in case $k\geq 2$ and the Heintze-Karcher inequality in case $k=1$,
imply that the flow \eqref{flow} starting from $M_{0}=M$ with $F = H_{k}/H_{k-1}$ preserves $W_{k}$ and increases $W_{k-1}$. As the flow converges to a geodesic sphere, at $t=\8$ the equality in \eqref{thm:QM-A} holds and hence the inequality is evident for $M=M_{0}$.
\end{proof}

Hence it remains to prove \autoref{thm:flow}, which is the main contribution of this paper.

As another byproduct of our main results we obtain a set of geometric inequalities for weighted quermassintegrals, based on an observation by Kwong/Wei \cite{KwongWei:10/2023}, who proved it in the convex class for the case $k=n$.

\begin{cor}
Let $M\sub \bbS^{n+1}_{+}$ be horo-convex. Then there holds
\eq{\label{thm:KW-A}(n+1)W_{k+1}(M)-\int_{M}\cos rH_{k}\geq \rho_{k}(W_{k}(M)),\q 1\leq k\leq n,}
where $\rho_{k}$ is a function which gives equality on geodesic spheres.
The equality case is precisely achieved on geodesic spheres.
\end{cor}

\pf{
In \cite[Lemma~5.5]{KwongWei:10/2023} it is computed that under the flow \eqref{flow} with $F = H_{k}/H_{k-1}$ the left hand side is decreasing, while it is known that $W_{k}(M)$ is preserved. Also the equality case is characterised. The convergence to a sphere completes the proof.
}

\subsection*{Some history}
Curvature flows have long been known to be powerful tools to prove geometric inequalities for hypersurfaces. The idea is that once a curvature flow is known to improve a geometric quantity, such as an isoperimetric ratio of any sort, and the limit of the flow and its value of the geometric quantity is known, a geometric inequality for the initial hypersurface follows. Natural candidates for such an endeavour are flows which preserve one quantity and decrease another one, such as the classical volume preserving mean curvature flow,
\eq{\del_{t}x = \br{\fint_{M_{t}}H - H}\nu,}
where $\fint_{M_{t}}H$ is the mean value integral of the mean curvature $H$. It preserves the enclosed volume, while decreasing the surface area. This flow was first studied by Huisken \cite{Huisken:/1987} for convex hypersurfaces of the Euclidean space and by Cabezas-Rivas/Miquel \cite{Cabezas-RivasMiquel:/2007} for horo-convex hypersurfaces of the hyperbolic space. Due to an unfavourable term in the evolution of the second fundamental form, the convexity of the hypersurface is not preserved for this flow in spherical space, and hence until today there are no results about this flow in the sphere. Similar flows where studied by Guan/Li \cite{GuanLi:08/2009} to prove the quermassintegral inequalities for not necessarily convex hypersurfaces of the Euclidean space,
\eq{W_{k}(M) \geq c_{n,k}W_{k-1}(M)^{\fr{n+1-k}{n+2-k}},}
where the $W_{k}$ are the natural counterparts of those from \eqref{Wk} in the Euclidean space. For this, they also employed a curvature flow containing a global term. Global terms such as the mean value of mean curvature make the analysis particularly hard, because some maximum principle arguments fail and there is no avoidance principle. Especially in non-Euclidean spaces these flow can lack smooth estimates. Hence about a decade ago Guan/Li \cite{GuanLi:/2015} invented a new type of curvature flow, which contains no global quantities, but still fixes a geometric quantity. Their prototype was the volume preserving mean curvature type flow
\eq{\del_{t}x = (\phi'-uH_{1})\nu,}
which also decreases area. The argument is similar to the one given in the proof of \autoref{thm:QM} above. The big advantage, compared to the standard volume preserving mean curvature flow is, that it does not require the convexity, but only starshapedness, and it works well in a large class of warped product spaces (see \cite{GuanLiWang:/2019}), including hyperbolic and spherical space.  
In particular, it is used to prove the isoperimetric inequality in a large class of warped product spaces of non-constant curvature. Recently,  Li and the first author  \cite{LIPan:/2025} extended this flow to more general Riemannian manifold admitting a non-trivial conformal vector field. The strength of this flow gave rise to the hope, that other versions of this, such as
\eq{\label{GL flow}\del_{t}x = \br{\phi'-u\fr{H_{k}}{H_{k-1}}}\nu}
or
\eq{\label{inverse GL flow}\del_{t}x = \br{\fr{\phi'H_{k-1}}{H_{k}}-u}\nu,}
which still share very good monotonicity properties towards the quermassintegral inequalities, could be employed to prove new such inequalities in space forms. In the literature, flows of this type are often called {\it locally constrained curvature flows}. The full conjecture by Guan/Li is, that these flows should converge for starshaped and $k$-convex initial hypersurfaces (i.e. $H_{i}>0$ for $i\leq k$).  Notably, even in the convex class the convergence of either \eqref{GL flow} and \eqref{inverse GL flow} is an open problem in the spherical and hyperbolic space. The full set of quermassintegral inequalities as in \autoref{thm:QM} is still open as well. In the hyperbolic space, they were proven for horo-convex hypersurfaces by Wang/Xia \cite{WangXia:07/2014} with a flow involving a global term and reproved using \eqref{inverse GL flow} by Hu/Li/Wei \cite{HuLiWei:04/2022}. Other special cases are known for the convex and positive sectional curvature class in the hyperbolic space, for example when $k=n$ or $k=1$. We do not give a precise account here, but mention \cite{AndrewsChenWei:/2021,AndrewsHuLi:/2020,HuLi:04/2019}. In the sphere, results are similarly scattered. For convex hypersurfaces, there were several special cases of \eqref{thm:QM-A} proved using varying methods. For example, Makowski and the second author \cite{MakowskiScheuer:11/2016} and Wei/Xiong \cite{WeiXiong:/2015} used purely expanding curvature flows to prove a subclass of inequalities. In \cite{ChenGuanLiScheuer:/2022} further special cases involving $W_{0}$ were treated using \eqref{GL flow} and also inequalities involving $W_{n}$ are possible, since in the sphere the flows \eqref{GL flow} and \eqref{inverse GL flow} are equivalent in the convex class due to duality. Until today, for $k\notin\{1,n\}$, neither of the flows \eqref{GL flow} or \eqref{inverse GL flow} is understood in any satisfactory way.

Despite the complicated behaviour of the flow \eqref{inverse GL flow} in the sphere, the quermassintegral inequalities for convex hypersurfaces are somewhat more complete than in the hyperbolic space. In \cite{ChenSun:03/2022}, Chen/Sun prove an inequality of the form
\eq{W_{k}(M)\geq \psi_{k}(W_{k-2}(M)),}
for all $k$.

 However, for the convex class in spherical space, \eqref{thm:QM-A} is still open. In addition, contrary to the hyperbolic space, until today no stronger convexity condition on a hypersurface of the sphere was known, such that \eqref{thm:QM-A} holds. In this paper, we provide such a condition and we call it horo-convexity. By this, we mirror the hyperbolic picture to the spherical one and the outcome is \autoref{thm:flow} and \autoref{thm:QM}.

\subsection*{Funding}
This work was funded by the German Research Foundation (Deutsche Forschungsgemeinschaft, DFG) through the project ``Curvature flows with local and nonlocal sources'', SCHE 1879/4-1.

\section{Horo-convexity and Evolution equations}

\subsection*{Horo-convexity}

In this section we motivate our definition of horo-convexity, which is clearly motivated geometrically by the following proposition.

\begin{prop}\label{prop:horo-spheres}
Let $M$ be a geodesic sphere, which lies in the northern hemisphere and touches the equator $\del\bbS^{n+1}_{+}$. Then on $M$ there holds
\eq{\label{horo}\phi'h = (1-u)g,}
where $\phi(r) = \sin(r)$, $u = \ip{\phi\del_{r}}{\nu}$ and $r$ is the radial distance to the north pole.
\end{prop}

\pf{
Let $\mathcal{O}$ be the north pole. Assume that $M$ is a geodesic sphere of radius $r_0$ centred at the point $p_0$. For any point $p\in M$, there is a unique geodesic triangle $\Delta(\mathcal{O},p_0,p)$ formed by geodesic arcs.
 Then, the spherical law of cosines gives
\eq{\label{consin}\cos d(\mathcal{O},p_0)=\cos d(\mathcal{O},p)\cos d(p,p_0)+\sin d(\mathcal{O},p)\sin d(p,p_0)\langle \nu(p),\partial_r\rangle,}
where  \(d(\cdot,\cdot)\) denotes the spherical distance, and \(\nu(p)\) is unit outer normal vector  of $M$ at the vertex \(p\). 
Note that \(d(p,p_0)=r_0\),  \(d(\mathcal{O},p_0)=r(p_0)\) and \(d(\mathcal{O},p)=r(p)\). 
Thus \eqref{consin} becomes
\begin{equation}\label{cosin-2}
	\cos r(p_0)
	= \cos r(p)\cos r_0 + \sin r_0\,u(p).
\end{equation}
Since \(M\) touches the equator \(\partial \mathbb{S}^{n+1}_{+}\), we have
\eq{
r(p_0)+r_0=\tfrac{\pi}{2}.
}
It follows that \(\cos r(p_0)=\sin r_0\), and hence \eqref{cosin-2} yields
\eq{\phi'(r(p))\frac{\cos r_0}{\sin r_0}+u(p)=1,}
which is \eqref{horo}.
}

From this property, it is natural to call spheres as in \autoref{prop:horo-spheres} the {\it horo-spheres} of $\bbS^{n+1}_{+}$ and the wording horo-convexity for a hypersurface is clearly motivated by a local one-sidedness.

\subsection*{Evolution Equations}
\begin{rem}[Short-time existence]
We start the flow \eqref{flow} from a horo-convex hypersurface in $\bbS^{n+1}_{+}$, in particular it is strictly convex and thus it has principal curvatures in $\Ga_{+}$. The curvature function $F$ gives rise to an elliptic operator, which is uniformly elliptic on compact subsets of $\Ga_{+}$. As $\phi'>0$, the requirements for short-time existence for a parabolic PDE are met and thus there exists a maximal time $T^{*}>0$ of smooth existence, compare \cite[Sec.~2.5]{Gerhardt:/2006} for details. This automatically means that $F$ is positive up to $T^{*}$. In the following, when we refer to the flow \eqref{flow}, we mean that time runs up to $T^{*}$, unless stated otherwise and we usually call $M_{0}$ the initial hypersurface.
\end{rem}
\begin{nota}
	Throughout this paper, $F$ is viewed as a function on the Weingarten operator $h^i_j=g^{ik}h_{jk}$, i.e $F=F(h^i_j)$, or equivalently as a function of $h$ and $g$, $F = F(h_{ij},g_{ij})$. We use the notation
	\eq{F^{i}_j=\frac{\partial F}{\partial h^i_j}=g_{jk}F^{ik} = g_{jk}\fr{\del F}{\del h_{ik}},\ F^{ij,k\ell}=\frac{\partial^2F}{\partial h_{ij}\partial h_{k\ell}}.}
We refer to \cite[Sec.~2.1]{Gerhardt:/2006} for details.
\end{nota}

For convenience, we define the function
\eq{\Phi(r)=\int_0^r\phi(s)ds=1-\cos r,}
then $\bar\n\Phi=\phi(r)\partial_r$ is a conformal Killing vector field on $\bbS^{n+1}$, where $\bar\n$ is the Levi-Civita connection of the round metric $\bar g$.

To shorten formulas and notation, in the following we use semi-colons to denote indices of covariant derivatives, for example
\eq{T_{k;ij} = (\n^{2}_{ij}T)_{k},}
where $\n$ is the Levi-Civita connection of the induced metric $g$ of a hypersurface $M\sub \bbS_{+}^{n+1}$. 

\begin{lemma}
Along the flow \eqref{flow} there hold the following evolution equations.
\eq{\label{ev-s}\del_{t}u&=  \fr{\phi'}{F^{2}}F^{ij}u_{;ij}+ \ip{\bar\n\Phi}{\n u}+\fr{\phi'}{F^{2}}(F^{ij}h_{ik}h^{k}_{j}-F^{2})u + \fr{\phi^{2}-u^{2}}{F}, }
\eq{\label{ev-thetaprime}\del_{t}\phi'&=\fr{\phi'}{F^{2}}F^{ij}\phi'_{;ij} + u^{2} - \fr{2\phi'u}{F} + \fr{\phi'^{2}}{F^{2}}F^{k}_{k}}
and
\eq{\label{ev-h}\del_{t}h^{i}_{j} &= \fr{\phi'}{F^{2}}F^{kl}h^{i}_{j;kl}+ \bar g(\bar\n\Phi,{x_{;}}^{k}){h^{i}_{j;k}} + \fr{1}{F^{2}}F_{;j}{\phi'_{;}}^{i} + \fr{1}{F^{2}}{F_{;}}^{i}\phi'_{;j}- \fr{2\phi'}{F^{3}}F_{;j}{F_{;}}^{i} \\
			&\hp{=}+ \fr{\phi'}{F^{2}}F^{kl,rs}{h_{kl;}}^{i}h_{rs;j}  + \fr{\phi'}{F^{2}}F^{kl}h^{r}_{k}h_{rl}h^{i}_{j} + \br{\phi'-\fr{\phi'}{F^{2}}F^{k}_{k}-\fr{u}{F}}h^{i}_{j} - \fr{2\phi'}{F}h^{i}_{k}h^{k}_{j} + \br{\fr{\phi'}{F}+u}\de^{i}_{j}.
} 
\end{lemma}

\pf{
According to \cite[p.~4722]{GuanLi:/2015} we have
\eq{\label{D2s}u_{;ij}&=\phi'h_{ij}-h_{ik}h^{k}_{j}u+\ip{\bar\n\Phi}{x_{;k}}h^{k}_{i;j}.}
Also there holds
\eq{\del_{t}u&=\phi'\br{\fr{\phi'}{F}-u}-\ip{\bar\n\Phi}{\n\br{\fr{\phi'}{F}-u}}=\fr{\phi'^{2}}{F}-u\phi'-\fr{1}{F}\ip{\bar\n\Phi}{\n\phi'}+\fr{\phi'}{F^{2}}\ip{\bar\n\Phi}{\n F}+\ip{\bar\n\Phi}{\n u},}
where we used the conformal Killing property of $\bar\n \Phi$ and the evolution of the outer normal $\nu$, cf. \cite[Lemma~2.3.2]{Gerhardt:/2006}.
Hence the first equation follows by combining these two relations with the Codazzi equation, the homogeneity of $F$ and
\eq{\ip{\bar\n\Phi}{\n\phi'} = -\phi^{2}g(\n r,\n r) = u^{2}-\phi^{2}.}

For the second equation we note that
\eq{\phi'(r)=\ip{x}{e_{n+2}},}
when $x$ is viewed as a codimension $2$ embedding into $\bbR^{n+2}$. Hence
\eq{\del_{t}\phi'=\br{\fr{\phi'}{F}-u}\ip{\nu}{e_{n+2}}}
and
\eq{\phi'_{;ij}=\ip{x_{;ij}}{e_{n+2}}=-h_{ij}\ip{\nu}{e_{n+2}}-g_{ij}\phi',}
where we have used the codimension two Gaussian formula
\eq{x_{;ij}=-h_{ij}\nu-g_{ij}x.}
Thus, using $\ip{\nu}{e_{n+2}}=-u$, we get
\eq{\del_{t}\phi'-\fr{\phi'}{F^{2}}F^{ij}\phi'_{;ij}&=\br{\fr{\phi'}{F}-u}\ip{\nu}{e_{n+2}}+\fr{\phi'}{F}\ip{\nu}{e_{n+2}}+\fr{\phi'}{F^{2}}F^{ij}g_{ij}\phi'=u^{2}-2u\fr{\phi'}{F}+\fr{\phi'}{F^{2}}F^{ij}g_{ij}\phi'.}

For the second fundamental form we proceed with \cite[Lemma~2.3.3]{Gerhardt:/2006}, and use the above formulas for $u_{;ij}$ and $\phi'_{;ij}$ together with some immediate cancellations to obtain
\eq{\del_{t}h^{i}_{j}&={\br{u-\fr{\phi'}{F}}_{;j}}^{i}+\br{u-\fr{\phi'}{F}}h^{i}_{k}h^{k}_{j}+\br{u-\fr{\phi'}{F}}\de^{i}_{j}\\
			&=\phi'h^{i}_{j} + \bar g(\bar\n\Phi,x_{;k}){h^{ik}}_{;j} + \fr{1}{F^{2}}F_{;j}{\phi'_{;}}^{i} + \fr{1}{F^{2}}{F_{;}}^{i}\phi'_{;j} + \fr{\phi'}{F^{2}}{F_{;j}}^{i} - \fr{2\phi'}{F^{3}}F_{;j}{F_{;}}^{i} - \fr{u}{F}h^{i}_{j} -\fr{\phi'}{F}h^{i}_{k}h^{k}_{j} + u\de^{i}_{j}.
} 
We have to expand the second derivatives of $F$. The standard Simon's identity is
\eq{F_{;ij} = F^{kl}h_{ij;kl} + F^{kl,rs}h_{kl;i}h_{rs;j} + F^{kl}h^{r}_{k}h_{rl}h_{ij} - F^{k}_{k}h_{ij} - Fh_{ik}h^{k}_{j} + Fg_{ij},}
see for example \cite[Equ.~(2.28)]{ChenGuanLiScheuer:/2022},
and plugging this in  gives the result.
}

Now we compute the evolution of the crucial tensor which defines the horo-convexity.

\begin{lemma}
The tensor 
\eq{S^{i}_{j} = \phi'h^{i}_{j} + (u-1)\de^{i}_{j}}
satisfies 
\eq{\del_{t}S^{i}_{j}
			&=\fr{\phi'}{F^{2}}F^{kl}S^{i}_{j;kl} - \fr{2}{F^{2}}F^{kl}\phi'_{;k}S^{i}_{j;l}+ \bar g(\bar\n\Phi,{x_{;}}^{k}){S^{i}_{j;k}} + \fr{\phi'}{F^{2}}F_{;j}{\phi'_{;}}^{i} + \fr{\phi'}{F^{2}}{F_{;}}^{i}\phi'_{;j}\\
			&\hp{=}- \fr{2\phi'^{2}}{F^{3}}F_{;j}{F_{;}}^{i}+ \fr{\phi'^{2}}{F^{2}}F^{kl,rs}{h_{kl;}}^{i}h_{rs;j}+ \fr{2}{F^{2}}F^{kl}\phi'_{;k}\phi'_{;l}h^{i}_{j}\\
			&\hp{=} + \br{\fr{\phi'}{F^{2}}F^{kl}h^{r}_{k}h_{rl}+\fr{1}{\phi'}+\fr{u-4}{F}}S^{i}_{j} - \fr{2}{F}S^{i}_{k}S^{k}_{j} + \br{\fr{\phi'}{F^{2}}F^{kl}h^{r}_{k}h_{rl}+ \fr{2}{F^{2}}F^{kl}\phi'_{;k}u_{;l} +\fr{1-u}{\phi'}+\fr{u-1}{F}} \de^{i}_{j}. 
			}
\end{lemma}

\pf{
By \eqref{ev-s}, \eqref{ev-thetaprime} and \eqref{ev-h},   we compute
\eq{\del_{t}S^{i}_{j} &= \del_{t}\phi' h^{i}_{j} + \phi'\del_{t}h^{i}_{j} + \del_{t}u \de^{i}_{j}\\
			&=\fr{\phi'}{F^{2}}F^{kl}S^{i}_{j;kl} - \fr{2\phi'}{F^{2}}F^{kl}\phi'_{;k}h^{i}_{j;l} + \br{ u^{2} - \fr{2\phi'u}{F} + \fr{\phi'^{2}}{F^{2}}F^{k}_{k}} h^{i}_{j}\\
			&\hp{=} + \phi'\bar g(\bar\n\Phi,{x_{;}}^{k}){h^{i}_{j;k}} + \fr{\phi'}{F^{2}}F_{;j}{\phi'_{;}}^{i} + \fr{\phi'}{F^{2}}{F_{;}}^{i}\phi'_{;j}- \fr{2\phi'^{2}}{F^{3}}F_{;j}{F_{;}}^{i} \\
			&\hp{=}+ \fr{\phi'^{2}}{F^{2}}F^{kl,rs}{h_{kl;}}^{i}h_{rs;j}  + \fr{\phi'^{2}}{F^{2}}F^{kl}h^{r}_{k}h_{rl}h^{i}_{j} + \br{\phi'^{2}-\fr{\phi'^{2}}{F^{2}}F^{k}_{k}-\fr{u\phi'}{F}}h^{i}_{j} - \fr{2\phi'^{2}}{F}h^{i}_{k}h^{k}_{j} + \br{\fr{\phi'^{2}}{F}+u\phi'}\de^{i}_{j}\\
			&\hp{=} + \br{ \bar g(\bar\n\Phi,\n u)+\fr{\phi'}{F^{2}}F^{rs}h_{rk}h^{k}_{s}u - u\phi' + \fr{\phi^{2}-u^{2}}{F}} \de^{i}_{j}. 
			}
We replace gradients of $h^{i}_{j}$ by those of $S^{i}_{j}$ as follows,

\eq{
- \fr{2\phi'}{F^{2}}F^{kl}\phi'_{;k}h^{i}_{j;l} = - \fr{2}{F^{2}}F^{kl}\phi'_{;k}S^{i}_{j;l} + \fr{2}{F^{2}}F^{kl}\phi'_{;k}\phi'_{;l}h^{i}_{j} + \fr{2}{F^{2}}F^{kl}\phi'_{;k}u_{;l}\de^{i}_{j}
}			
and
\eq{\phi'\bar g(\bar\n\Phi,{x_{;}}^{k}){h^{i}_{j;k}}&=\bar g(\bar\n\Phi,{x_{;}}^{k}){S^{i}_{j;k}} - \bar g(\bar\n\Phi,\n \phi')h^{i}_{j} - \bar g(\bar\n \Phi,\n u)\de^{i}_{j}\\
				&=\bar g(\bar\n\Phi,{x_{;}}^{k}){S^{i}_{j;k}} + \phi^{2}g(\n r,\n r)^{2}h^{i}_{j} - \bar g(\bar\n \Phi,\n u)\de^{i}_{j}\\
				&=\bar g(\bar\n\Phi,{x_{;}}^{k}){S^{i}_{j;k}} + (\phi^{2}-u^{2})h^{i}_{j} - \bar g(\bar\n \Phi,\n u)\de^{i}_{j}.}

We plug this in and already rearrange a bit, while using $\phi^{2} + \phi'^{2}=1$,

\eq{\del_{t}S^{i}_{j}
			&=\fr{\phi'}{F^{2}}F^{kl}S^{i}_{j;kl} - \fr{2}{F^{2}}F^{kl}\phi'_{;k}S^{i}_{j;l}+ \bar g(\bar\n\Phi,{x_{;}}^{k}){S^{i}_{j;k}} \\
			&\hp{=}+ \fr{\phi'}{F^{2}}F_{;j}{\phi'_{;}}^{i} + \fr{\phi'}{F^{2}}{F_{;}}^{i}\phi'_{;j}- \fr{2\phi'^{2}}{F^{3}}F_{;j}{F_{;}}^{i}+ \fr{\phi'^{2}}{F^{2}}F^{kl,rs}{h_{kl;}}^{i}h_{rs;j} + \fr{2}{F^{2}}F^{kl}\phi'_{;k}\phi'_{;l}h^{i}_{j} + \fr{2}{F^{2}}F^{kl}\phi'_{;k}u_{;l}\de^{i}_{j}\\
			&\hp{=} + \br{\fr{\phi'^{2}}{F^{2}}F^{kl}h^{r}_{k}h_{rl}+1-\fr{3u\phi'}{F}}h^{i}_{j} - \fr{2\phi'^{2}}{F}h^{i}_{k}h^{k}_{j} + \br{\fr{\phi'}{F^{2}}F^{rs}h_{rk}h^{k}_{s}u + \fr{1-u^{2}}{F}} \de^{i}_{j}. 
			}

For the next step, we replace $\phi'h^{i}_{j} = S^{i}_{j} + (1-u)\de^{i}_{j}$ except the term involving $\phi'_{;k}$,
hence we obtain

\eq{\del_{t}S^{i}_{j}
			&=\fr{\phi'}{F^{2}}F^{kl}S^{i}_{j;kl} - \fr{2}{F^{2}}F^{kl}\phi'_{;k}S^{i}_{j;l}+ \bar g(\bar\n\Phi,{x_{;}}^{k}){S^{i}_{j;k}} \\
			&\hp{=}+ \fr{\phi'}{F^{2}}F_{;j}{\phi'_{;}}^{i} + \fr{\phi'}{F^{2}}{F_{;}}^{i}\phi'_{;j}- \fr{2\phi'^{2}}{F^{3}}F_{;j}{F_{;}}^{i}+ \fr{\phi'^{2}}{F^{2}}F^{kl,rs}{h_{kl;}}^{i}h_{rs;j} + \fr{2}{F^{2}}F^{kl}\phi'_{;k}\phi'_{;l}h^{i}_{j} + \fr{2}{F^{2}}F^{kl}\phi'_{;k}u_{;l}\de^{i}_{j}\\
			&\hp{=} + \br{\fr{\phi'}{F^{2}}F^{kl}h_{rk}h^{r}_{l}+\fr{1}{\phi'}-\fr{3u}{F}}(S^{i}_{j} + (1-u)\de^{i}_{j}) - \fr{2}{F}(S^{i}_{k} + (1-u)\de^{i}_{k})(S^{k}_{j} + (1-u)\de^{k}_{j})\\
			&\hp{=} + \br{\fr{\phi'}{F^{2}}F^{kl}h_{rk}h^{r}_{l}u + \fr{1-u^{2}}{F}} \de^{i}_{j}\\
			&=\fr{\phi'}{F^{2}}F^{kl}S^{i}_{j;kl} - \fr{2}{F^{2}}F^{kl}\phi'_{;k}S^{i}_{j;l}+ \bar g(\bar\n\Phi,{x_{;}}^{k}){S^{i}_{j;k}} + \fr{\phi'}{F^{2}}F_{;j}{\phi'_{;}}^{i} + \fr{\phi'}{F^{2}}{F_{;}}^{i}\phi'_{;j}\\
			&\hp{=}- \fr{2\phi'^{2}}{F^{3}}F_{;j}{F_{;}}^{i}+ \fr{\phi'^{2}}{F^{2}}F^{kl,rs}{h_{kl;}}^{i}h_{rs;j} + \fr{2}{F^{2}}F^{kl}\phi'_{;k}\phi'_{;l}h^{i}_{j}+ \fr{2}{F^{2}}F^{kl}\phi'_{;k}u_{;l}\de^{i}_{j}\\
			&\hp{=} + \br{\fr{\phi'}{F^{2}}F^{kl}h^{r}_{k}h_{rl}+\fr{1}{\phi'}+\fr{u-4}{F}}S^{i}_{j} - \fr{2}{F}S^{i}_{k}S^{k}_{j} + \br{\fr{\phi'}{F^{2}}F^{kl}h^{r}_{k}h_{rl}+\fr{1-u}{\phi'}+\fr{u-1}{F}} \de^{i}_{j}. 
			}
}

\section{Preservation of horo-convexity}

In this section,  we prove the main key to our results, namely the preservation of horo-convexity. Also note that we do not a priori assume the hypersurface to be starshaped around the north pole, but we recall that the initial hypersurface lies in the northern hemisphere $\bbS^{n+1}_{+}$.
We first prove that the maximum distance to north pole is decreasing.

\begin{prop}\label{lemma:phi' bound}
Let $M_{0}\sub \bbS^{n+1}_{+}$ be a horo-convex initial hypersurface, then the maximal distance of the flow hypersurfaces to the north pole is strictly decreasing, as long as it  is not the constant flow of centred spheres.
\end{prop}

\pf{
We use the evolution equation \eqref{ev-thetaprime} of $\phi'(r) = \cos r$, which can be rewritten as
\eq{\del_{t}\phi' = \fr{\phi'}{F^{2}}F^{ij}\phi'_{;ij} + \br{u-\fr{\phi'}{F}}^{2} + \fr{\phi'^{2}}{F^{2}}(F^{k}_{k}-1)\geq \fr{\phi'}{F^{2}}F^{ij}\phi'_{;ij},}
 because $F^{k}_{k}\geq 1$ due to the concavity of $F$. The result follows from the strong maximum principle, since if the maximum would not strictly decrease, $\phi'(r)$ would have to be constant.
}

Next, to complete the proof, we need the following lemma, which was proved by Brendle, Choi, and Daskalopoulos \cite[Lemma 5]{BrendleChoiDaskalopoulos:/2017} and by Choi, Kim, and Lee \cite[Lemma 4.1]{ChoiKimLee:2024}.

\begin{lemma}\label{lemma: viscosity}
	Let $D$ be the multiplicity of the smallest eigenvalue at a point $\xi_0$ on $M_{t_0}$ for $t_0 > 0$ so that 
	\eq{
	\sigma_1 = \cdots = \sigma_{D}< \sigma_{D+1} \le \cdots \le \sigma_n,
	}
	where $\sigma_1, \cdots, \sigma_n$ are the  eigenvalues of a tensor $\mathcal{P}$. Suppose $\eta$ is a smooth function defined on 
	\eq{
	\mathcal{M} = \bigcup_{0 < t \le t_0} M_t \times \{t\}
	}
	such that $\eta \le \sigma_1$ on $\mathcal{M}$ and $\eta= \sigma_1$ at $ (t_0,\xi_0)$. Then, at the point $(t_0,\xi_0)$ with coordinates satisfying 
	$g_{ij} = \delta_{ij}$ and $\mathcal{P}_{ij} = \sigma_i \delta_{ij}$, we have 
	\eq{\mathcal{P}_{kl;i} = \partial_i \eta \, \delta_{kl}, \quad \text{for} \quad 1 \le k,l \le D,}
	 and
	 \eq{ \eta_{;ii} \le \mathcal{P}^{1}_{1;ii} - 
		\sum_{l> D} \frac{2 (\mathcal{P}^{1}_{l;i})^2}{\sigma_{l}- \sigma_1}, 
		\quad
		\partial_t \eta\ge \partial_t \mathcal{P}^1_{1}.
}
\end{lemma}

\begin{thm}\label{thm:pres h-conv}
Along the flow \eqref{flow}, if the tensor $S$ has only non-negative eigenvalues at time zero, then this property is preserved.
\end{thm}

\pf{
Let $\si_{m} = \phi'\ka_{m} + (u-1)$ be the eigenvalues of $S$ and $D$ be the multiplicity of $\si_{1}$. 
We will show that the set of times in $[0,T^{*})$, where $\si_{1}\geq 0$ is an open set in $[0,T^{*})$. By continuity of $\si_{1}$, it is also closed and this will complete the proof. Without loss of generality we prove that $\si_{1}\geq0$ on the time interval $0\leq t\leq\delta$ where $\delta$ is small compared to a constant $C$ depending on $\max \vert h\vert$ and $\max\fr{1}{\kappa_1}$. We define
\eq{\tilde{S}_{ij}=S_{ij}+\varepsilon(\delta+t)g_{ij},}
and  we claim that $\tilde{S}_{ij}>0$ on $0\leq t\leq\delta$ for every $\varepsilon$. Then letting $\varepsilon\to0$ will finish the proof.

If not, there will be a first time $t_0$ with $0<t_0\leq\delta$ where the smallest eigenvalue of $\tilde{S}$ equals zero at some point $\xi_0\in\mathcal{M}_{t_0}$. Therefore, 
$\eta\equiv0$ is a smooth lower support of $\tilde{\sigma}_1$ at the point $(t_{0},\xi_{0})$, i.e. $\eta= \tilde{\si}_1(t_{0},\xi_{0})$ and  $\eta\leq \tilde{\si}_1$ on $\mathcal{M} = \bigcup_{0 < t \le t_0} M_t \times \{t\}$, where $\tilde{\si}_i=\si_i+\varepsilon(\delta+t)$ are the eigenvalues of the tensor $\tilde{S}$.

Then by Lemma \ref{lemma: viscosity}, at $(t_{0},\xi_{0})$ we obtain
\eq{\label{gradient S}\tilde{S}^{i}_{j;k} =0, \q \text{for} \quad 1 \le i,j  \le D,}
\eq{\label{maximum S}\del_{t}\tilde{S}^{1}_{1} \leq 0, \q \tilde{S}^{1}_{1;kk} - 2\sum_{l>D}\fr{(\tilde{S}^{1}_{l;k})^{2}}{\tilde{\si}_{l}-\tilde{\si}_{1}}\geq0.}

We have to plug  the evolution of $S^{1}_{1}$ into   \eqref{maximum S}, and then, we can ignore all gradients of $S^{1}_{1}$ by \eqref{gradient S}. Note that $\tilde{S}^1_1=\tilde{\si}_1=0$ at $(t_{0},\xi_{0})$, which means $\sigma_1(t_0,\xi_0)=-\varepsilon(\delta+t_0)\geq-2\varepsilon\delta$. Next, we keep $\varepsilon,\delta\leq1$ and compute in an orthonormal frame at $(t_{0},\xi_{0})$, 
\eq{\label{pf:pres h-conv 3}-\varepsilon&\geq   2\fr{\phi'}{F^{2}}F_{;1}{\phi'_{;1}} - \fr{2\phi'^{2}}{F^{3}}(F_{;1})^{2}+ \fr{\phi'^{2}}{F^{2}}F^{kl,rs}{h_{kl;1}}h_{rs;1}+ \fr{2\phi'^{2}}{F^{2}}F^{kk}\sum_{l>D}\fr{(h_{1l;k})^{2}}{\ka_{l}-\ka_{1}}\\
			&\hp{=}+ \fr{2}{F^{2}}F^{kl}\phi'_{;k}u_{;l}+ \fr{2}{F^{2}}F^{kl}\phi'_{;k}\phi'_{;l}\ka_{1} + \br{\fr{\phi'}{F^{2}}F^{kl}h^{r}_{k}h_{rl} +\fr{1-u}{\phi'}+\fr{u-1}{F}}-C\varepsilon\delta,   }
where we used, that for all $l>D$ and all $k$,
\eq{S_{1l;k} = \phi'_{;k}h_{1l} + \phi'h_{1l;k} + u_{;k}\de_{1l} = \phi'h_{1l;k},}
and $C$ is a constant depending on $\max \vert h\vert$ and $\max\fr{1}{\kappa_1}$.
We rewrite the following terms:
\eq{\label{pf:pres h-conv 1}&\fr{2\phi'^{2}}{F^{2}}F^{kk}\sum_{l>D}\fr{(h_{1l;k})^{2}}{\ka_{l}-\ka_{1}}+ \fr{2}{F^{2}}F^{kl}\phi'_{;k}u_{;l}+ \fr{2}{F^{2}}F^{kl}\phi'_{;k}\phi'_{;l}\ka_{1}\\
=~&\fr{2\phi'^{2}}{F^{2}}\sum_{k\leq D,l>D}F^{kk}\fr{(h_{1l;k})^{2}}{\ka_{l}-\ka_{1}}+\fr{2\phi'^{2}}{F^{2}}\sum_{k>D,l>D}F^{kk}\fr{(h_{1l;k})^{2}}{\ka_{l}-\ka_{1}}- \fr{2\phi^{2}}{F^{2}}F^{kl}r_{;k}r_{;l}\ka_{l}+ \fr{2\phi^{2}}{F^{2}}F^{kl}r_{;k}r_{;l}\ka_{1}\\
=~&\fr{2\phi'^{2}}{F^{2}}\sum_{k>D,l>D}F^{kk}\fr{(h_{1l;k})^{2}}{\ka_{l}-\ka_{1}} -\fr{2\phi^{2}}{F^{2}}(F^{kk}-F^{11})r_{;k}^{2}(\ka_{k}-\ka_{1}),}
where we used, for $k\leq D$,
\eq{\label{pf:pres h-conv 5}\phi'h_{1l;k} = \phi'h_{1k;l} = S_{1k;l} - \phi'_{;l}h_{1k} - u_{;l}\de_{1k} =\phi r_{;l}(\ka_{1}-\ka_{l})\delta_{1k}.}

We use a well known relation for the second derivatives of a curvature function, see \cite[Lemma 2.1.14]{Gerhardt:/2006} for details,
\eq{\label{pf:pres h-conv 2}F^{kl,rs}{h_{kl;1}}h_{rs;1} &= \fr{\del^{2} F}{\del\ka_{k}\del\ka_{l}}h_{kk;1}h_{ll;1} + 2\sum_{k>l}\fr{F^{kk}-F^{ll}}{\ka_{k}-\ka_{l}}(h_{kl;1})^{2}\\
						&= \sum_{k,l> D}\fr{\del^{2} F}{\del\ka_{k}\del\ka_{l}}h_{kk;1}h_{ll;1} + 2\sum_{l> D, k>l}\fr{F^{kk}-F^{ll}}{\ka_{k}-\ka_{l}}(h_{kl;1})^{2} + 2\sum_{l\leq D,k>l}\fr{F^{kk}-F^{ll}}{\ka_{k}-\ka_{l}}(h_{1l;k})^{2} \\
						&=\sum_{k,l > D} F^{kl,rs}h_{kl;1}h_{rs;1} + 2\sum_{l\leq D,k>l}\fr{F^{kk}-F^{ll}}{\ka_{k}-\ka_{l}}(h_{1l;k})^{2}.}
Next, from the inverse concavity we may use \cite[Prop.~4.3]{Scheuer:06/2018} to conclude
\eq{\label{pf:pres h-conv 4}F^{kl,rs}{h_{kl;1}}h_{rs;1} 
	&\geq \fr{2}{F}(F_{;1})^{2} - 2\sum_{k,l> D}F^{kr}b^{l s}h_{kl;1}h_{rs;1}+2\sum_{l\leq D,k>l}\fr{F^{kk}-F^{ll}}{\ka_{k}-\ka_{l}}(h_{1l;k})^{2}, }
	where $b$ is the inverse of $h$, and we used $h_{kk,1}=0$ for all $1\leq k\leq D$ by $\eqref{pf:pres h-conv 5}$.
We combine \eqref{pf:pres h-conv 1} and \eqref{pf:pres h-conv 4} to obtain
\eq{
&- \fr{2\phi'^{2}}{F^{3}}(F_{;1})^{2}+ \fr{\phi'^{2}}{F^{2}}F^{kl,rs}{h_{kl;1}}h_{rs;1}+ \fr{2\phi'^{2}}{F^{2}}F^{kk}\sum_{l>D}\fr{(h_{1l;k})^{2}}{\ka_{l}-\ka_{1}}+ \fr{2}{F^{2}}F^{kl}\phi'_{;k}u_{;l}+ \fr{2}{F^{2}}F^{kl}\phi'_{;k}\phi'_{;l}\ka_{1}\\
\geq~&- \fr{2\phi'^{2}}{F^{2}}\sum_{k,l> D}\fr{F^{kk}}{\ka_{l}}(h_{kl;1})^{2}+\fr{2\phi'^{2}}{F^{2}}\sum_{l\leq D,k>l}\fr{F^{kk}-F^{ll}}{\ka_{k}-\ka_{l}}(h_{1l;k})^{2}\\
	&+\fr{2\phi'^{2}}{F^{2}}\sum_{k,l>D}F^{kk}\fr{(h_{1l;k})^{2}}{\ka_{l}-\ka_{1}} -\fr{2\phi^{2}}{F^{2}}(F^{kk}-F^{11})r_{;k}^{2}(\ka_{k}-\ka_{1})\\
	=~& \fr{2\phi'^{2}}{F^{2}}\sum_{k,l> D}\fr{F^{kk}\ka_{1}}{\ka_{l}(\ka_{l}-\ka_{1})}(h_{kl;1})^{2}\\
	\geq~&\fr{2\phi'^{2}}{F^{2}}\sum_{k> D}\fr{F^{kk}\ka_{1}}{\ka_{k}(\ka_{k}-\ka_{1})}(h_{kk;1})^{2},
}
where in the equality we used \eqref{pf:pres h-conv 5} again.

Now we estimate
\eq{(F_{;1})^{2}\leq (\sum_{k>D}F^{kk}\abs{h_{kk;1}})^{2}\leq \sum_{k>D}\fr{F^{kk}\ka_{1}}{\ka_{k}(\ka_{k}-\ka_{1})}(h_{kk;1})^{2}\sum_{k>D}F^{kk}\fr{\ka_{k}(\ka_{k}-\ka_{1})}{\ka_{1}}}
and plug everything back into \eqref{pf:pres h-conv 3},

\eq{-\varepsilon&\geq   2\fr{\phi'}{F^{2}}F_{;1}{\phi'_{;1}}+ \fr{2\phi'^{2}}{F^{2}}\sum_{k> D}\fr{F^{kk}\ka_{1}}{\ka_{k}(\ka_{k}-\ka_{1})}(h_{kk;1})^{2} + \br{\fr{\phi'}{F^{2}}F^{kl}h^{r}_{k}h_{rl} +\fr{1-u}{\phi'}+\fr{u-1}{F}}-C\varepsilon\delta\\
			&\geq   2\fr{\phi'}{F^{2}}F_{;1}{\phi'_{;1}}+ \fr{2\phi'^{2}\ka_{1}}{F^{2}}\fr{(F_{;1})^{2}}{F^{kk}\ka_{k}^{2} - F\ka_{1}} + \br{\fr{\phi'}{F^{2}}F^{kl}h^{r}_{k}h_{rl} +\fr{1-u}{\phi'}+\fr{u-1}{F}}-C\varepsilon\delta\\
			&\geq    -\fr{F^{kk}\ka_{k}^{2}-F\ka_{1}}{2F^{2}\ka_{1}}(\phi'_{;1})^{2} + \br{\fr{\phi'}{F^{2}}F^{kl}h^{r}_{k}h_{rl} +\fr{1-u}{\phi'}+\fr{u-1}{F}}-C\varepsilon\delta.}
			
Finally, we use $\si_1 = \phi'\ka_{1} + u-1\geq-2\varepsilon\delta$ and discard $(1-u)/
\phi'>0$ to obtain
\eq{\label{pf:pres h-conv 6}-\varepsilon&\geq -\fr{F^{kk}\ka_{k}^{2}-F\ka_{1}}{2F^{2}\ka_{1}}(\phi'_{;1})^{2} + \fr{\phi'}{F^{2}}F^{kl}h^{r}_{k}h_{rl}-\fr{\phi'\ka_{1}}{F}-C\varepsilon\delta\geq\fr{(F^{kk}\ka_{k}^{2}-F\ka_{1})}{2F^{2}\ka_{1}}\br{2(1-u) - (\phi'_{;1})^{2}}-C\varepsilon\delta.}
Since 
\eq{2-2u-(\phi'_{;1})^{2}=2-2u-\phi^{2}r_{;1}^{2}\geq2-2u - \phi^{2}+u^{2} = (1-u)^{2}+1-\phi^{2}\geq 0,}
this yields $\varepsilon\leq C\varepsilon\delta$, which gives a contradiction when $C\delta<1$.
The proof is complete.
}

\section{Long-time existence and starshapedness}

In this section we prove $T^{*}=\8$ and that from some time onwards the north pole will remain inside the convex bodies enclosed by the flow hypersurfaces.
We begin with a simple observation.

\begin{lemma}
Let the initial hypersurface $M_{0}\sub \bbS^{n+1}_{+}$ be horo-convex, then the smallest principal curvature and the curvature function are uniformly bounded from below by a positive constant depending on the initial data.
\end{lemma}

\pf{
The preservation of horo-convexity implies
\eq{h\geq \fr{1-u}{\phi'}g\geq (1-\max_{M_{t}}\phi)g\geq (1-\max_{M_{0}}\phi)g,}
as $\max \phi$ is strictly decreasing due to \autoref{lemma:phi' bound}.
Hence, the smallest principal curvature is bounded below.
Moreover, since  $F$ is  homogeneous 1 and strictly monotone increasing in each argument, there holds
\[F(\kappa_1,\cdots,\kappa_n)\geq F(\kappa_1,\cdots,\kappa_1)=\kappa_1,\]
which provides a lower bound for $F$.
}

We obtain a time-dependent bound on the second fundamental form as follows.

\begin{lemma}
Let the initial hypersurface $M_{0}\sub \bbS^{n+1}_{+}$ be horo-convex, then along the flow \eqref{flow} the second fundamental form is bounded by a constant depending only on initial data and on $T^{*}$.
\end{lemma}

\pf{
Due to a well known trick (see \cite[Lemma 4.5]{Gerhardt:/2014} for details), it suffices to estimate the evolution of $h^{n}_{n}$. From \eqref{ev-h}, the lower $F$ bound, Young's inequality and the concavity of $F$ we obtain at spatial maxima,
\eq{\del_{t}h^{n}_{n}&\leq \fr{2}{F^{2}}F_{;n}\phi_{;n}' - \fr{2\phi'}{F^{3}}(F_{;n})^{2} + \fr{\phi'}{F^{2}}F^{kl}h^{r}_{k}h_{rl}\ka_{n} - \fr{2\phi'}{F}\ka_{n}^{2} + c(\ka_{n}+1)\leq  \fr{\phi'}{F}\ka_{n}^{2} - \fr{2\phi'}{F}\ka_{n}^{2} + c(\ka_{n}+1).  }
The result follows from the maximum principle. Note that we can not use the second order term to obtain a time-independent bound, since we do not have the upper bound on $F$ yet.
}

\begin{cor}
Let the initial hypersurface $M_{0}\sub \bbS^{n+1}_{+}$ be horo-convex, then the flow exists for all times and for any given time $t$ it satisfies uniform $C^{\8}$-estimates up to $t$.
\end{cor}

\pf{
Suppose $T^{*}<\8$. Then up to $T^{*}$, all curvature quantities are uniformly bounded from above and below. Hence the convex flow hypersurfaces $M_{t}$ satisfy a uniform inradius bound from below. As the flow moves with finite velocity, there exists a point $p\in \bbS^{n+1}_{+}$, such that near $T^{*}$ the point $p$ is in the convex region enclosed by $M_{t}$ and there holds
\eq{\min_{M_{t}}\dist(p,\cdot)\geq c>0}
for some constant $c$.
Hence, all $M_{t}$ are starshaped around $p$ with uniform estimates on all quantities. Using the corresponding graph functions, the flow can be extended beyond $T^{*}$ in a standard way, compare for example \cite[Sec.~2.5]{Gerhardt:/2006}.
}

Now that we know the long-time existence, we can prove that the flow becomes starshaped around the north pole after some waiting time.

\begin{lemma}\label{lemma:starshaped}
Let the initial hypersurface $M_{0}\sub \bbS^{n+1}_{+}$ be horo-convex, then the flow becomes starshaped around the north pole after some waiting time and then remains a starshaped flow in the northern hemisphere. 
\end{lemma}

\pf{
We argue with certain properties of the related flow
\eq{\label{pf:starshaped 1}\del_{t}y = \br{\fr{\phi'}{H_{1}}-u}\nu,}
which preserves the surface area and increases the enclosed volume due to the Heintze-Karcher inequality
\eq{\int_{M_{t}}\fr{\phi'}{H_{1}}\geq \int_{M_{t}}u,}
which holds for all mean convex hypersurfaces of the northern hemisphere. Equality holds if and only if $M_{t}$ is a sphere. Now start the flow $(N_{t})$ solving \eqref{pf:starshaped 1} from any geodesic sphere inside $\bbS^{n+1}_{+}$. Hence, equality holds in the isoperimetric inequality for $N_{0}$, and therefore each subsequent $N_{t}$ must also be a sphere; otherwise, its enclosed volume would increase by the Heintze-Karcher inequality while its surface area remained constant. Hence the flow $(N_{t})$ satisfies
\eq{\del_{t}y = \br{\fr{\phi'}{H_{1}}-u}\nu = \br{\fr{\phi'}{F}-u}\nu,}
because $H_{1} = F$ for geodesic spheres. Therefore, $(N_{t})$ also satisfies \eqref{flow} and by uniqueness we have $M_{t} = N_{t}$. Consequently,  for \eqref{flow} with an initial sphere $S_{r_0}(p_0)$, the flow hypersurfaces are spheres of the same radius.
They converge to a sphere centred at the north pole, due to \autoref{lemma:phi' bound} and the fact that stationary points of the flow must be centred spheres.

Now let $M_{0}$ be a horo-convex initial hypersurface and $(M_{t})_{t>0}$ be the flow \eqref{flow} emerging from it. Pick an inball, the boundary of which moves to a sphere centred at the north pole by the previous argument. From the avoidance principle, after some time the north pole must be in the region bounded by $M_{t}$. 
}

\section{Estimates for horo-convex and starshaped flows}

Now that starshapedness and preservation of horo-convexity is settled for our flow, the remaining estimates for completion of the convergence are fairly straightforward. All we have to do is refine our previous estimates for the second fundamental form.

\begin{cor}
Let the initial hypersurface $M_{0}\sub \bbS^{n+1}_{+}$ be horo-convex.
Then the largest principal curvature, and hence the full second fundamental form, is bounded from above during the evolution \eqref{flow}.
\end{cor}

\pf{
The hypersurface $M_{0}$ is contained in the open set $\bbS^{n+1}_{+}$. Hence, due to \autoref{lemma:phi' bound}, there holds $\phi'\geq \min_{M_{0}}\phi'>0$.
By a standard trick, wlog assume $\ka_{n} = h^{n}_{n}$ and consider the evolution of
\eq{w = \log h^{n}_{n} - \log u.}
Then, at a maximum of $w$ and with the help of the concavity of $F$,
\eq{\del_{t}w - \fr{\phi'}{F^{2}}F^{ij}w_{;ij}&\leq \fr{2\phi'}{F\ka_{n}}(\log F)_{;n}{(\log\phi')_{;n}} - \fr{2\phi'}{F\ka_{n}}(\log F)_{;n}^{2}\\
			&\hp{=}+2\phi'  - \fr{\phi'}{F^{2}}F^{k}_{k}-\fr{u}{F} - \fr{2\phi'}{F}\ka_{n} + \br{\fr{\phi'}{F}+u}\ka_{n}^{-1} - \fr{\phi^{2}-u^{2}}{uF}\\
			&\leq \fr{\phi'}{F\ka_{n}}(\log\phi')_{;n}^{2} - \fr{\phi^{2}-u^{2}}{uF}\\
			&\leq\fr{\phi^2-u^2}{\phi'uF\ka_{n}}(1-\phi'\kappa_n), }
where in the second inequality we used the facts $F\leq F(\kappa_n,\cdots,\kappa_n)=\kappa_n$ and $F^k_k\geq1$.
Then, at a maximum of $\omega$, $\del_{t}w\leq0$ if $\ka_{n}$ is large enough. Hence, $\ka_{n}$ is bounded. 
}

\begin{cor}
Let the initial hypersurface $M_{0}\sub \bbS^{n+1}_{+}$ be horo-convex. Then the flow \eqref{flow} starting from $M_{0}$ exists for all times and converges smoothly to a geodesic sphere centred at the north pole.
\end{cor}

\pf{
Written as graphs around the north pole, the radial function of the flow hypersurfaces $M_{t}$ satisfies a fully nonlinear elliptic equation of the form
\eq{\del_{t}r = G(D^{2}r,Dr,r,\cdot),}
where the precise form of $G$ can explicitly be computed from $F$ the geometric quantities of $M_{t}$. Our geometric bounds yield uniform estimates up to $C^{2}$, and they also yield the uniform parabolicity of the curvature function $F$ along the evolution. Standard parabolic regularity theory gives parabolic H\"older estimates of the second derivatives and then Schauder theory applied to the linearisation of the operator gives estimates to arbitrary order.  Due to these uniform estimates, a diagonal argument and Arzela-Ascoli, the sequence of flows 
\eq{x_k(t,\cdot)=x(t+k,\cdot)}
subsequentially converges to a limit flow $x_{\infty}$ with corresponding radial function $\phi'_{\infty}$.  By Prop \ref{lemma:phi' bound}, the  quantity $\min_{M}\phi'(t,x(t,\cdot))$ is increasing along the flow, and therefore converges to a constant  $c_0$. Hence, there holds
\eq{\min_{M}\phi'_{\infty}(t,x_{\infty}(t,\cdot))=\lim_{k\to\infty}\min_{M}\phi'(t,x(t+k,\cdot))=c_0.}
By the strong maximum principle, this implies $\phi'_{\infty}\equiv c_0$. Consequently, the flow converges to a sphere centred at the north pole.
}

\providecommand{\bysame}{\leavevmode\hbox to3em{\hrulefill}\thinspace}
\providecommand{\href}[2]{#2}

\end{document}